\documentclass[preprint,12pt]{elsarticle}
\usepackage[utf8]{inputenc}
\usepackage{amsmath}
\usepackage{amssymb}
\usepackage{amsthm}
\usepackage{xcolor}

\theoremstyle{plain}
\newtheorem{theorem}{Theorem}
\newtheorem{lemma}[theorem]{Lemma}

\theoremstyle{definition}
\newtheorem{definition}[theorem]{Definition}
\newtheorem{example}[theorem]{Example}
\newtheorem{remark}[theorem]{Remark}

\begin{document}

\begin{frontmatter}

\title{Nordhaus--Gaddum relations for degrees and connectivities under the $\delta$-complement}
\author[aff1]{Supakorn Srisawat}
\ead{supakorn.swt@gmail.com}

\author[aff1,aff2]{Panupong Vichitkunakorn\corref{cor}} 
\ead{panupong.v@psu.ac.th}

\cortext[cor]{Corresponding author}

\affiliation[aff1]{organization={Division of Computational Science, Faculty of Science},
            addressline={Prince of Songkla University}, 
            city={Hat Yai, Songkhla},
            postcode={90110}, 
            country={Thailand}}

\affiliation[aff2]{organization={Research Center in Mathematics and Statistics with Applications},
            addressline={Prince of Songkla University}, 
            city={Hat Yai, Songkhla},
            postcode={90110}, 
            country={Thailand}}

\begin{abstract}
The $\delta$-complement $G_\delta$ of a graph $G$ complements adjacency within each degree class and preserves adjacency between distinct degree classes.
We establish sharp lower and upper additive and product Nordhaus--Gaddum bounds for the minimum degree, maximum degree, vertex connectivity, and edge connectivity of $G$ and $G_\delta$.
We give explicit constructions attaining all stated sharp bounds.
\end{abstract}


\end{frontmatter}

\section{Introduction}
Nordhaus and Gaddum~\cite{noga} initiated the study of inequalities relating a graph invariant on a graph $G$ and on its complement $\overline G$.
These inequalities are now called \emph{Nordhaus--Gaddum relations}.
They have been studied for minimum degree~\cite{alavi}, maximum degree~\cite{xu2}, diameter, girth, and circumference~\cite{xu}, domination parameters~\cite{japa}, and many more.
A survey of Nordhaus--Gaddum relations is given in~\cite{aouchiche}.

Pai et al.~\cite{pai} introduced the $\delta$-complement $G_\delta$, which complements the subgraph induced by each degree class of $G$ and leaves adjacency between distinct degree classes unchanged.
This operation is related to subgraph complementation~\cite{kaminski2009}, which complements a specified induced subgraph.
It is also related to Seidel switching~\cite{van1966}, which reverses adjacency between a specified vertex set and its complement.

Nordhaus--Gaddum relations for the chromatic numbers of $G$ and $G_\delta$ were studied by Vichitkunakorn et al.~\cite{panupong}.
Recent work has considered domination parameters of $G$ and $G_\delta$~\cite{tangjai2026,chen2026}.
Here we establish sharp lower and upper additive and product Nordhaus--Gaddum bounds for the minimum degree, maximum degree, vertex connectivity, and edge connectivity of $G$ and $G_\delta$.
We also give explicit constructions attaining all stated sharp bounds.

\section{Preliminaries}\label{bg}
All graphs in this paper are finite, simple, and undirected.
For a graph $G$, write $n=|V(G)|$, and let $\deg_G(v)$ denote the degree of $v$.
The minimum and maximum degrees are denoted by $\delta(G)$ and $\Delta(G)$.
The vertex connectivity $\kappa(G)$ is the minimum number of vertices whose deletion disconnects $G$ or leaves a single vertex, and the edge connectivity $\lambda(G)$ is the minimum number of edges whose deletion disconnects $G$.
We set $\kappa(G)=\lambda(G)=0$ when $G$ is disconnected or has order one.
The disjoint union and join of graphs $G$ and $H$ are denoted by $G+H$ and $G\vee H$, respectively.
As usual, $mG$ denotes the disjoint union of $m$ copies of $G$.
For a degree $d$ occurring in $G$, let
\[
V_d(G)=\{v\in V(G):\deg_G(v)=d\}.
\]

\begin{definition}[Pai et al.~\cite{pai}]\label{defdelcom}
The $\delta$-complement of $G=(V,E)$ is the graph $G_\delta=(V,E_\delta)$ in which distinct vertices $u$ and $v$ are adjacent if and only if either
\begin{enumerate}
\item $\deg_G(u)=\deg_G(v)$ and $uv\notin E$, or
\item $\deg_G(u)\neq\deg_G(v)$ and $uv\in E$.
\end{enumerate}
\end{definition}

Thus $G_\delta[V_d(G)]=\overline{G[V_d(G)]}$ for every degree class, whereas adjacency between distinct degree classes is unchanged.
We shall also use the following compatibility with ordinary complementation.

\begin{lemma}[\cite{pai}]\label{comp_deltacomp_lemma}
For every graph $G$, we have \(\overline{G_\delta}\cong(\overline G)_\delta.\)
\end{lemma}

\section{Degree and Connectivity Relations}\label{ng}

We begin with the minimum degree.

\begin{definition}\label{def:sSa}

For $p\geq1$ and $1\leq i \leq a$, define
\[
s_p(i) = \begin{cases} i+p-2, &\text{if $i$ and $p$ are even},\\ i+p-1, &\text{otherwise.}\end{cases}
\]
Then define its partial sum for $a\geq 1$ as
\[
S_p(a)= \sum_{i=1}^a s_p(i) = 
\begin{cases}
\displaystyle\sum_{i=1}^{a}(i+p-1)-\left\lfloor a/2\right\rfloor,
&\text{if $p$ is even},\\[2mm]
\displaystyle\sum_{i=1}^{a}(i+p-1),
&\text{if $p$ is odd}.
\end{cases}
\]
For $n,p\geq1$, define
\[
a_n(p)=\min\{a\geq1:S_p(a)\geq n\}.
\]

\end{definition}

\begin{lemma}\label{Spformulas}
For $p,n\geq1$,
\begin{equation}\label{anpdef}
a_n(p)=
\begin{cases}
\left\lceil\sqrt{(p-1)^2+2n-1}-(p-1)\right\rceil,
&p\text{ even},\\[2mm]
\left\lceil\sqrt{(p-\tfrac12)^2+2n}-(p-\tfrac12)\right\rceil,
&p\text{ odd}.
\end{cases}
\end{equation}
\end{lemma}

\begin{proof}
If $p$ is odd, then
\[
S_p(a)=\sum_{i=1}^{a}(i+p-1)=\frac{a(a+2p-1)}2.
\]
Solving the resulting quadratic inequality gives
\[
S_p(a)\geq n
\quad\Longleftrightarrow\quad
a\geq
\sqrt{\left(p-\frac12\right)^2+2n}
-\left(p-\frac12\right).
\]
The least integral $a$ is the ceiling of the expression on the right.

If $p$ is even, then
\[
\begin{aligned}
S_p(a)
&=\frac{a(a+2p-1)}2-\left\lfloor\frac a2\right\rfloor\\
&=\frac{a(a+2p-2)}2
+\frac{a-2\lfloor a/2\rfloor}{2}
=\left\lceil\frac{a(a+2p-2)}2\right\rceil.
\end{aligned}
\]
Indeed, the second term on the middle line is zero when $a$ is even and is $1/2$ when $a$ is odd.
The first term is respectively an integer or a half-integer.
Therefore the last equality holds.
Moreover,
\[
\begin{split}
S_p(a)\geq n
&\quad\Longleftrightarrow\quad
\frac{a(a+2p-2)}2>n-1\\
&\quad\Longleftrightarrow\quad
a(a+2p-2)>2n-2\\
&\quad\Longleftrightarrow\quad
a(a+2p-2)\geq2n-1\\
&\quad\Longleftrightarrow\quad
(a+p-1)^2\geq(p-1)^2+2n-1.
\end{split}
\]
The third line follows because $a(a+2p-2)$ is an integer.
Taking the positive square root shows that the least integral $a$ is
\[
\left\lceil\sqrt{(p-1)^2+2n-1}-(p-1)\right\rceil,
\]
which proves \eqref{anpdef}.
\end{proof}

\begin{lemma}\label{mindegbound}
Let $G$ be a graph of order $n$ such that $\delta(G_\delta)=n-p$, where $1\leq p\leq n$.
Then
\[
\delta(G)\leq n-a_n(p).
\]
\end{lemma}

\begin{proof}
Put $a=n-\delta(G)$.
Since every vertex has degree at least $\delta(G)=n-a$, each nonempty degree class has index $i\in\{1,\ldots,a\}$.
Fix $i\in\{1,\ldots,a\}$.
If $V_{n-i}(G)$ is nonempty, choose $v\in V_{n-i}(G)$ and put
\[
b(v)=\deg_{G[V_{n-i}(G)]}(v).
\]
Every neighbor of $v$ in $V_{n-i}(G)$ is a nonneighbor of $v$ in $G_\delta$, so $b(v)\leq p-1$.
The vertex $v$ has $n-i-b(v)$ neighbours outside $V_{n-i}(G)$, whereas only $n-|V_{n-i}(G)|$ vertices lie outside that degree class.
Hence
\[
|V_{n-i}(G)|\leq i+b(v)\leq i+p-1.
\]
Suppose that $p$ and $i$ are both even and that equality holds in the final inequality.
Since the displayed inequality holds for every vertex of $V_{n-i}(G)$, equality forces $b(v)=p-1$ for every such vertex.
Thus $G[V_{n-i}(G)]$ is $(p-1)$-regular of order $i+p-1$, which is impossible because both its order and its degree are odd.
Therefore, 
\begin{equation}\label{eq:sizeV}
|V_{n-i}(G)| \leq \begin{cases} i+p-2, &\text{if $i$ and $p$ are even},\\ i+p-1, &\text{otherwise.}\end{cases}
\end{equation}
Definition~\ref{def:sSa} gives $|V_{n-i}(G)| \leq s_p(i)$.
Thus
\[
n=\sum_{i=1}^{a}|V_{n-i}(G)|
\leq\sum_{i=1}^{a}s_p(i)=S_p(a).
\]
Hence $a\geq a_n(p)$.
\end{proof}

\begin{theorem}\label{ubsum}
Let $G$ be a graph of order $n\geq1$, and let $k$ be the integer satisfying
\[
\binom{k}{2}<n\leq\binom{k+1}{2}.
\]
Then
\begin{equation}\label{mindegsum}
0\leq\delta(G)+\delta(G_\delta)\leq2n-k-1.
\end{equation}
For $n\neq2$, a graph with exactly one isolated vertex attains the lower bound.
For every $n$, there is a set $I=\{i_1,\dots,i_t\}\subseteq\{1,\ldots,k\}$ such that $k\in I$ and $\sum_{i\in I}i=n$, and the complete multipartite graph $K_{i_1, \dots,i_t}$ attains the upper bound.
\end{theorem}

\begin{proof}
The lower bound is immediate.
If $G$ has exactly one isolated vertex, that vertex forms a singleton degree class and remains isolated in $G_\delta$, so the lower bound is attained whenever $n\neq2$.

For the upper bound, put $a=n-\delta(G)$ and $p=n-\delta(G_\delta)$.
By Lemma~\ref{mindegbound}, $a\geq a_n(p)$.
Since $S_p$ is increasing, the definition of $a_n(p)$ gives $n\leq S_p(a)$.
Therefore
\[
n\leq S_p(a)\leq\frac{a(a+2p-1)}2
=\binom{a+p}{2}-\binom p2
\leq\binom{a+p}{2}.
\]
Since $n>\binom k2$, it follows that $a+p\geq k+1$.
Therefore
\[
\delta(G)+\delta(G_\delta)
=2n-(a+p)\leq2n-k-1.
\]

To obtain an extremal graph, observe that
\[
0\leq n-k\leq\binom{k}{2}
\]
and that the subset sums of $\{1,\ldots,k-1\}$ comprise every integer from $0$ to $\binom{k}{2}$.
Choose a subset with sum $n-k$ and adjoin $k$ to form $I$.
The complete multipartite graph with part orders in $I$ has minimum degree $n-k$.
Its part orders are distinct, so its partite sets are precisely its degree classes and its $\delta$-complement is $K_n$.
The sum in \eqref{mindegsum} is therefore $2n-k-1$.
\end{proof}

We next determine the maximum product of the two minimum degrees.

\begin{lemma}[\cite{erdosgallai}]\label{regulargraphexistence}
	Let $n\geq1$ and $d\geq0$.
	A $d$-regular graph of order $n$ exists if and only if $0\leq d\leq n-1$ and $nd$ is even.
\end{lemma}

The lemma gives the following construction.
For $d=0$, take the empty graph on vertex set $\mathbb Z_n$.
If $d>0$ is even, join each $i\in\mathbb Z_n$ to $i\pm j\pmod n$ for $1\leq j\leq d/2$.
If $d$ is odd, then $n$ is even.
Join each $i\in\mathbb Z_n$ to $i\pm j\pmod n$ for $1\leq j\leq(d-1)/2$ and to $i+n/2\pmod n$.
In each case, every vertex has degree $d$.

\begin{theorem}\label{ubprod}
For every graph $G$ of order $n\geq2$,
\begin{equation}\label{exactDformula}
0\leq\delta(G)\delta(G_\delta)
\leq\max_{1\leq p\leq n}(n-p)\bigl(n-a_n(p)\bigr).
\end{equation}
The lower bound is attained by $K_n$.
The upper bound is attained for every $n$ by a join of regular graphs.
\end{theorem}

\begin{proof}
The lower bound follows from the nonnegativity of minimum degrees.
Moreover, $(K_n)_\delta=\overline{K_n}$, so $K_n$ attains the lower bound.

For an arbitrary graph $G$ of order $n$, set
\[
p=n-\delta(G_\delta).
\]
By Lemma~\ref{mindegbound},
\[
\delta(G)\delta(G_\delta)
\leq (n-p)\bigl(n-a_n(p)\bigr).
\]
Taking the maximum over $p\in\{1,\ldots,n\}$ proves the upper bound in \eqref{exactDformula}.

It remains to construct a graph attaining that bound.
For $n=2$, the right-hand side of \eqref{exactDformula} is zero, and $K_2$ is an extremal regular graph.
Hence assume that $n\geq3$.
Choose $p\in\{1,\ldots,n\}$ and $a=a_n(p)$ such that $(n-p)(n-a)$ is the maximum in \eqref{exactDformula}.
The maximizing term is positive, since $(K_{1,n-1})_\delta=K_n$ and therefore $\delta(K_{1,n-1})\delta((K_{1,n-1})_\delta)=n-1$.
Consequently, $n-a>0$ and $n-p>0$.

Our goal is to construct a join of regular graphs
\[
H=\bigvee_{i=1}^{a} R_i,
\]
where, whenever $w_i>0$, $R_i$ is a $d_i$-regular graph of order $w_i=i+d_i$ and $0\leq d_i\leq p-1$.
The numbers $w_1,\ldots,w_a$ are nonnegative integers satisfying $\sum_{i=1}^{a}w_i=n$ and $w_a>0$.
An empty block $R_i$ is omitted from the join.

Let $v$ be a vertex in a nonempty $R_i$.
These properties give $\deg_H(v)=n-i$.
Therefore $\delta(H)=n-a$.
The $\delta$-complement replaces $R_i$ by $\overline{R_i}$ and gives $\deg_{H_\delta}(v)=n-1-d_i\geq n-p$. 
Therefore $\delta(H_\delta)\geq n-p$.
The upper bound already proved will force $\delta(H_\delta)=n-p$.

For $1\leq i\leq a$, let
\[
\mathcal W_i=\{0\}\cup
\{i+d:0\leq d\leq p-1,\ (i+d)d\text{ is even}\}.
\]
The positive elements of $\mathcal W_i$ begin with $i$.
If $i$ is odd, they are all the integers from $i$ to $i+p-1$.
If $i$ is even, they are $i,i+2,i+4,\ldots$, ending with $i+p-2$ when $p$ is even and with $i+p-1$ when $p$ is odd.
So the gap between two consecutive positive elements is at most two.
Let $m_i=\max\mathcal W_i$ and $M_i=\sum_{j=1}^{i}m_j$.
From \eqref{eq:sizeV}, for $1\leq i\leq a$, we have
\[
m_i = \begin{cases} i+p-2, &\text{if $i$ and $p$ are even},\\ i+p-1, &\text{otherwise.}\end{cases}
\]
For $2\leq i\leq a$, we have $m_i= s_p(i)=S_p(i)-S_p(i-1)$, and $M_i=S_p(i)$ for every $1\leq i\leq a$.

To show that $n$ can be written as $\sum_{i=1}^{a}w_i$ with $w_i\in\mathcal W_i$, we will prove by induction on $j$ the stronger statement that every integer in $[0,M_j]$ can be written as $\sum_{i=1}^{j}w_i$ with $w_i\in\mathcal W_i$.

For $j=1$, we have $\mathcal W_1=\{0,1,\ldots,p\}$, so the assertion holds.

Now let $j\geq2$, and suppose that the first $j-1$ sets generate every integer in $[0,M_{j-1}]$.
Write elements in $\mathcal W_j$ as $x_0=0$ and $j =  x_1 < x_2 < \cdots < x_r = m_j$.
By the induction hypothesis, choosing $w_j=x_t$ produces every integer in $x_t+[0,M_{j-1}]$.
We need to show that
\[
\bigcup_{t=0}^r \left(x_t + [0, M_{j-1}] \right) = [0, M_j].
\]
It is enough to show that 
\begin{equation}\label{eq:xt_xt+1}
x_{t+1} \leq x_t + M_{j-1} + 1 \quad\text{for all } 0 \leq t < r.
\end{equation}
Since $j-1\in\mathcal W_{j-1}$, we have $M_{j-1}\geq m_{j-1}\geq j-1.$
So \eqref{eq:xt_xt+1} holds for $t=0$.
Moreover, $M_{j-1}\geq m_1=p$.
For $t\geq1$, we have $x_{t+1}-x_t\leq2\leq p+1\leq M_{j-1}+1$.
Hence \eqref{eq:xt_xt+1} also holds for $t\geq1$.
This completes the induction.

Choose $w_i\in\mathcal W_i$ with $\sum_{i=1}^{a}w_i=n$.
For every $i$ with $w_i>0$, write $w_i=i+d_i$.
Then $0\leq d_i\leq w_i-1$ and $w_id_i$ is even, so Lemma~\ref{regulargraphexistence} supplies a $d_i$-regular graph $R_i$ of order $w_i$.
Lastly, to show $w_a > 0$, we write $n = \sum_{i=1}^{a}w_i$.
If $w_a=0$, then
\[
n=\sum_{i=1}^{a-1}w_i\leq M_{a-1}=S_p(a-1)<n,
\]
which is impossible.
Thus $w_a>0$.

The preceding degree calculation now gives $\delta(H)=n-a$ and $\delta(H_\delta)=n-p$, so $H$ attains the upper bound in \eqref{exactDformula}.
\end{proof}

\begin{example}
Let $n=11$.
The values of the expression on the right-hand side of \eqref{exactDformula} are as follows.
\[
\begin{array}{c|ccccccccccc}
p&1&2&3&4&5&6&7&8&9&10&11\\ \hline
(11-p)(11-a_{11}(p))&60&63&64&56&54&45&36&27&18&9&0
\end{array}
\]
The largest value is $64$, attained when $p=3$, for which $a_{11}(3)=3$.
For this choice of $p$, the admissible sets are
\[
\mathcal W_1=\{0,1,2,3\},
\qquad
\mathcal W_2=\{0,2,4\},
\qquad
\mathcal W_3=\{0,3,4,5\}.
\]
One admissible representation is $11=2+4+5=(1+1)+(2+2)+(3+2)$
where $w_1=2, w_2=4, w_3=5$ and $d_1=1, d_2=2, d_3=2$.
One example of regular graphs is $R_1=K_2$, $R_2=C_4$, and $R_3=C_5$.
Consequently,
\[
H=K_2\vee C_4\vee C_5.
\]
Moreover,
\[
H_\delta
=\overline{K_2}\vee\overline{C_4}\vee\overline{C_5}
\cong 2K_1\vee2K_2\vee C_5.
\]
Hence $\delta(H)=8$ and $\delta(H_\delta)=8$.
\[\delta(H)\delta(H_\delta)= 64 = (11-3)(11-a_{11}(3)).\]
\end{example}

We now consider the maximum degree.

\begin{theorem}\label{maxdegsum}
Let $G$ be a graph of order $n\geq1$, and let
\[
\binom{k}{2}<n\leq\binom{k+1}{2}.
\]
Then
\begin{equation}\label{maxdegreesum}
k-1\leq\Delta(G)+\Delta(G_\delta)\leq2n-2.
\end{equation}
There exists a set $I=\{i_1,\ldots,i_t\}\subseteq\{1,\ldots,k\}$ such that $k\in I$ and $\sum_{i\in I}i=n$, and $K_{i_1}+\cdots+K_{i_t}$ attains the lower bound.
The upper bound is attained by $K_1$ when $n=1$ and by $K_{1,n-1}$ when $n\geq3$.
\end{theorem}

\begin{proof}
For every graph $H$ of order $n$,
\[
\Delta(H)=n-1-\delta(\overline H).
\]
Lemma~\ref{comp_deltacomp_lemma} therefore gives
\[
\begin{aligned}
\Delta(G)+\Delta(G_\delta)
&=2n-2-\delta(\overline G)-\delta(\overline{G_\delta})\\
&=2n-2-\delta(\overline G)-\delta((\overline G)_\delta).
\end{aligned}
\]
Applying Theorem~\ref{ubsum} to $\overline G$ proves \eqref{maxdegreesum}.

For the lower equality, let $I$ be a set supplied by Theorem~\ref{ubsum} and write $I=\{i_1,\ldots,i_t\}$.
The components of
\[
G=K_{i_1}+\cdots+K_{i_t}
\]
have distinct orders and hence form distinct degree classes.
It follows that $G_\delta=nK_1$, while $\Delta(G)=k-1$.
For $n\geq3$, the two degree classes of $K_{1,n-1}$ are its center and its leaves, and direct application of the definition gives $(K_{1,n-1})_\delta=K_n$.
Thus both maximum degrees are $n-1$.
The case $n=1$ is immediate.
\end{proof}

\begin{theorem}\label{maxdegreeproduct}
For every graph $G$ of order $n\geq1$,
\[
0\leq\Delta(G)\Delta(G_\delta)\leq(n-1)^2.
\]
The lower bound is attained by $K_n$.
The upper bound is attained by $K_1$ when $n=1$ and by $K_{1,n-1}$ when $n\geq3$.
\end{theorem}

\begin{proof}
Since
\[
0\leq\Delta(G),\Delta(G_\delta)\leq n-1,
\]
both inequalities follow.
Since $(K_n)_\delta=nK_1$, the complete graph attains the lower bound.
The upper-bound constructions are the same as in Theorem~\ref{maxdegsum}.
\end{proof}

\begin{remark}\label{ordertworemark}
For either graph of order two,
\[
\delta(G)+\delta(G_\delta)
=\Delta(G)+\Delta(G_\delta)=1
\]
and
\[
\delta(G)\delta(G_\delta)
=\Delta(G)\Delta(G_\delta)=0.
\]
Thus the lower additive bounds and the upper maximum-degree bounds have the stated order-two exception.
\end{remark}

We next consider additive relations for vertex and edge connectivity.

\begin{theorem}\label{connectivitysum}
Let $G$ be a graph of order $n\geq1$, and let
\[
\binom{k}{2}<n\leq\binom{k+1}{2}.
\]
Then
\[
0\leq\kappa(G)+\kappa(G_\delta)\leq2n-k-1
\]
and
\[
0\leq\lambda(G)+\lambda(G_\delta)\leq2n-k-1.
\]
For $n\neq2$, a graph with exactly one isolated vertex attains both lower bounds.
There is a set $I=\{i_1,\dots,i_t\}\subseteq\{1,\ldots,k\}$ such that $k\in I$ and $\sum_{i\in I}i=n$, and the complete multipartite graph $K_{i_1, \dots,i_t}$ attains both upper bounds.
\end{theorem}

\begin{proof}
The inequalities follow from Theorem~\ref{ubsum} and
\[
\kappa(H)\leq\lambda(H)\leq\delta(H).
\]
If $G$ has exactly one isolated vertex, then both $G$ and $G_\delta$ are disconnected, which gives the two lower equalities.

For the upper equalities, let $G$ be the complete multipartite graph from Theorem~\ref{ubsum}.
Its largest part has order $k$.
Deleting the $n-k$ vertices outside that part disconnects $G$ or leaves one vertex, while deleting fewer vertices leaves vertices in at least two parts and hence a connected graph.
Therefore $\kappa(G)=n-k$.
Also $\delta(G)=n-k$, so
\[
n-k=\kappa(G)\leq\lambda(G)\leq\delta(G)=n-k.
\]
The part orders are distinct, and hence $G_\delta=K_n$.
It follows that
\[
\kappa(G_\delta)=\lambda(G_\delta)=n-1,
\]
which proves both upper equalities.
\end{proof}

\begin{remark}
For either graph of order two,
\[
\kappa(G)+\kappa(G_\delta)
=\lambda(G)+\lambda(G_\delta)=1.
\]
\end{remark}

We now determine the maximum edge-connectivity product.

\begin{lemma}\label{denseedgelemma}
If a graph $H$ of order $n\geq1$ satisfies $\delta(H)\geq\lfloor n/2\rfloor$, then $\lambda(H)=\delta(H)$.
\end{lemma}

\begin{proof}
The assertion is immediate for $n=1$.
For $n\geq2$, the degree condition makes $H$ connected.
Let $S$ be the smaller side of an edge cut $[S,\overline{S}]$, and put $s=|S|$.
Then $1\leq s\leq\lfloor n/2\rfloor$, and
\[
|[S,\overline{S}]|\geq s\delta(H)-s(s-1)
=\delta(H)+(s-1)(\delta(H)-s)
\geq\delta(H).
\]
Thus every edge cut has at least $\delta(H)$ edges.
Together with $\lambda(H)\leq\delta(H)$, this proves the lemma.
\end{proof}

\begin{theorem}\label{exactedgeproduct}
For every graph $G$ of order $n\geq1$,
\begin{equation}\label{exactLformula}
0\leq\lambda(G)\lambda(G_\delta)
\leq\max_{1\leq p\leq n}(n-p)\bigl(n-a_n(p)\bigr).
\end{equation}
If $n=1$, both bounds are zero and are attained by $K_1$.
For $n\geq2$, the lower bound is attained by $K_n$.
For $n=2,\ldots,7$, the upper bound is attained respectively by
\[
2K_1,\quad K_{1,2},\quad K_2\vee\overline{K_2},\quad
K_{2,3},\quad K_{1,2,3},\quad K_{1,1,2,3}.
\]
For $n\geq8$, a join of regular graphs from the proof of Theorem~\ref{ubprod} attains the upper bound.
\end{theorem}

\begin{proof}
If $n=1$, then $G=K_1$, $G_\delta=K_1$, and both sides of \eqref{exactLformula} are zero.
Hence suppose that $n\geq2$.
The lower bound follows from the nonnegativity of edge connectivities.
Since $(K_n)_\delta=\overline{K_n}$, the graph $K_n$ attains the lower bound.

The upper bound follows from $\lambda(H)\leq\delta(H)$ and Theorem~\ref{ubprod}.
Suppose first that $n\geq8$, and let $k$ satisfy
\[
\binom{k}{2}<n\leq\binom{k+1}{2}.
\]
The complete multipartite construction in Theorem~\ref{ubsum} gives
\[
\max_{|V(H)|=n}\delta(H)\delta(H_\delta)\geq(n-k)(n-1).
\]
Let $H$ be a join of regular graphs from Theorem~\ref{ubprod}, and write
\[
x=\delta(H),
\qquad
y=\delta(H_\delta).
\]
If $x<n-k$, then
\[
xy\leq(n-k-1)(n-1)<(n-k)(n-1),
\]
This is a contradiction.
The same argument applies to $y$.
For $8\leq n\leq10$, one has $k=4$ and $n-k\geq\lceil n/2\rceil$.
For $k\geq5$, the inequality $n>\binom{k}{2}$ implies $n\geq2k$, and again $n-k\geq\lceil n/2\rceil$.
Hence
\[
x,y\geq n-k\geq\lceil n/2\rceil.
\]
Lemma~\ref{denseedgelemma}, applied to the extremal graph and its $\delta$-complement, gives $\lambda=\delta$ for both graphs.
Hence the upper bound in \eqref{exactLformula} is attained for $n\geq8$.
For orders two through seven, direct evaluation of $\delta(H)\delta(H_\delta)$ gives $0,2,4,8,15,20$, which attain the bound.
\end{proof}

Vertex connectivity requires additional restrictions.

\begin{definition}\label{def:sSak}
For $p,a\geq1$ and $1\leq i\leq a$, define
\begin{equation}\label{eq:sk}
s_{p,a}^\kappa(i)=
\begin{cases}
i+p-3=p-1,&\text{if $p\geq3$ is odd and $i=2$},\\
i+p-2=a,&\text{if $p=2$, $a\geq3$ is odd, and $i=a$},\\
i+p-2,&\text{if $p$ and $i$ are even},\\
i+p-1,&\text{otherwise}.
\end{cases}
\end{equation}
Then define its partial sum for $a\geq1$ as
\[
S_p^\kappa(a)=\sum_{i=1}^{a}s_{p,a}^\kappa(i)=
\begin{cases}
S_p(a)-1,&\text{if $p=2$ and $a\geq3$ is odd},\\
S_p(a)-2,&\text{if $p\geq3$ is odd and $a\geq2$},\\
S_p(a),&\text{otherwise}.
\end{cases}
\]
For $n,p\geq1$, define
\[
a_n^\kappa(p)=\min\{a\geq1:S_p^\kappa(a)\geq n\}.
\]
\end{definition}

\begin{lemma}\label{vertexconnectivitybound}
Let $G$ be a graph of order $n\geq1$ such that $\kappa(G_\delta)=n-p$, where $1\leq p\leq n$.
Then
\[
\kappa(G)\leq n-a_n^\kappa(p).
\]
\end{lemma}

\begin{proof}
Put $a=n-\kappa(G)$.
Since
\[
\delta(G)\geq\kappa(G)=n-a,
\]
every nonempty degree class is $V_{n-i}(G)$ for some $i\in\{1,\ldots,a\}$.
Consider $v\in V_{n-i}(G)$.
Let $b_v$ be the number of neighbors of $v$ in $V_{n-i}(G)$, and let $q_v$ be the number of nonneighbors of $v$ outside $V_{n-i}(G)$.
Counting nonneighbors of $v$ in the two graphs gives
\begin{equation}\label{vertexclassequations}
i-1=|V_{n-i}(G)|-1-b_v+q_v,
\qquad
n-1-\deg_{G_\delta}(v)=q_v+b_v.
\end{equation}
The second equation gives $q_v + b_v \leq p-1$ since $\delta(G_\delta) \geq \kappa(G_\delta) = n-p$.
The first equation gives 
\begin{equation}\label{eq:Vip}
	|V_{n-i}(G)|=i+b_v-q_v\leq i+b_v+q_v\leq i+p-1.
\end{equation}
Suppose that $p$ and $i$ are both even and that equality holds in $|V_{n-i}(G)|\leq i+p-1$.
For every $v\in V_{n-i}(G)$, both inequalities in \eqref{eq:Vip} must then be equalities.
The first equality gives $q_v=0$, and the second gives $b_v=p-1$.
Since $v$ is arbitrary, we see that $G[V_{n-i}(G)]$ is $(p-1)$-regular of odd order $i+p-1$, which is impossible.
Hence $|V_{n-i}(G)|\leq i+p-2$ when $p$ and $i$ are both even.

Definition~\ref{def:sSa} gives $|V_{n-i}(G)|\leq s_p(i)$ for every $1\leq i\leq a$.

The following arguments account for the two additional cases in \eqref{eq:sk}.

First suppose that $p\geq3$ is odd and $a\geq2$.
For $i=2$, the bound improves from $s_p(2)=p+1$ to $p-1=s_{p,a}^\kappa(2)$.
If $|V_{n-2}(G)|=p+1$, then \eqref{vertexclassequations} forces $q_v=0$ and $b_v=p-1$ for every $v\in V_{n-2}(G)$.
Hence $G_\delta[V_{n-2}(G)]$ is $1$-regular of order $p+1\ge4$.
It is therefore disconnected.
Then $V(G)\setminus V_{n-2}(G)$ is a vertex cut of $G_\delta$ of size $n-(p+1)=\kappa(G_\delta)-1$, a contradiction.
If $|V_{n-2}(G)|=p$, then \eqref{vertexclassequations} forces $q_v=0$ and $b_v=p-2$, which would give a regular graph of odd order $p$ and odd degree $p-2$.
Thus $|V_{n-2}(G)|\leq p-1=s_{p,a}^\kappa(2)$.

Now suppose that $p=2$ and $a\geq3$ is odd.
For $i=a$, the bound improves from $s_2(a)=a+1$ to $a=s_{2,a}^\kappa(a)$.
If $|V_{n-a}(G)|=a+1$, then \eqref{vertexclassequations} forces $q_v=0$ and $b_v=1$ for every $v\in V_{n-a}(G)$.
Thus $G[V_{n-a}(G)]$ is $1$-regular of order $a+1\ge4$ and is disconnected.
The set $V(G)\setminus V_{n-a}(G)$ is then a vertex cut of $G$ size $n-(a+1)=\kappa(G)-1$, a contradiction.
Hence $|V_{n-a}(G)|\leq a=s_{2,a}^\kappa(a)$.

Consequently, $|V_{n-i}(G)|\leq s_{p,a}^\kappa(i)$ for every $1\leq i\leq a$.
Therefore
\[
n=\sum_{i=1}^{a}|V_{n-i}(G)|
\leq\sum_{i=1}^{a}s_{p,a}^\kappa(i)
=S_p^\kappa(a).
\]
Thus $a\geq a_n^\kappa(p)$, and hence $\kappa(G)=n-a\leq n-a_n^\kappa(p)$.
\end{proof}

\begin{definition}
Let $s\geq3$ be an integer.
The circulant graph $C_s(a_1,\ldots,a_k)$ has vertex set $\mathbb Z_s$, where $0<a_1<\cdots<a_k\leq\lfloor s/2\rfloor$, and vertices $x$ and $y$ are adjacent when $x-y\equiv\pm a_j\pmod s$ for some $j$.
The integers $a_1,\ldots,a_k$ are called the jumps of the circulant.
\end{definition}

\noindent The following criterion identifies a class of maximally connected circulant graphs.

\begin{lemma}[\cite{boeschtindell}, Theorem~1]\label{consecutivejumpconnectivity}
Let $s\geq3$ and $1\leq a\leq b\leq\lfloor s/2\rfloor$ be integers.
If $C_s(a,a+1,\ldots,b)$ is connected, then
\[
\kappa\bigl(C_s(a,a+1,\ldots,b)\bigr)
=\delta\bigl(C_s(a,a+1,\ldots,b)\bigr).
\]
\end{lemma}

\noindent We apply this criterion to the complements of powers of cycles.

\begin{lemma}\label{cyclepowerconnectivity}
Let $s\geq5$ and $1\leq q\leq \frac{s-3}2$ be integers.
Let $C_s^q$ denote the $q$-th power of the cycle $C_s$.
\[
\kappa(C_s^q)=2q,
\qquad
\kappa(\overline{C_s^q})=s-1-2q.
\]
\end{lemma}
\begin{proof}
The first equality is stated explicitly in \cite[Lemma~3.1]{liliu2024}.
For the second equality, put $G=\overline{C_s^q}=C_s(q+1,\ldots,\lfloor s/2\rfloor)$.
If this jump sequence has at least two terms, then it contains consecutive integers and hence has greatest common divisor one with $s$.
If it has one term, then $s$ is odd, $q=(s-3)/2$, and the jump $(s-1)/2$ is relatively prime to $s$.
Thus $G$ is connected.
Lemma~\ref{consecutivejumpconnectivity} gives $\kappa(G)=\delta(G)=s-1-2q$.
\end{proof}

\noindent We now use the preceding bounds and connectivity facts to determine the maximum vertex-connectivity product.

\begin{theorem}\label{vertexproduct}
For every graph $G$ of order $n\geq1$,
\begin{equation}\label{exactKformula}
0\leq\kappa(G)\kappa(G_\delta)
\leq
\max_{1\leq p\leq n}(n-p)\bigl(n-a_n^\kappa(p)\bigr).
\end{equation}
The lower bound is attained by $K_n$, and the upper bound is attained for every order by a join of regular graphs.
\end{theorem}

\begin{proof}
For an arbitrary graph $G$, put $p=n-\kappa(G_\delta)$.
Lemma~\ref{vertexconnectivitybound} gives $\kappa(G)\leq n-a_n^\kappa(p)$, and multiplication by $\kappa(G_\delta)=n-p$ proves the upper bound in \eqref{exactKformula}.
The lower bound follows from the nonnegativity of vertex connectivity, and $K_n$ attains it because $(K_n)_\delta=\overline{K_n}$.

It remains to attain the bound.
Choose $p\in\{1,\ldots,n\}$ at which the maximum in \eqref{exactKformula} is attained, and put $a=a_n^\kappa(p)$.
For $1\leq i\leq a$, let
\[
\mathcal W_i=\{0\}\cup
\left\{i+d:\begin{array}{l}
0\leq d\leq p-1,\quad (i+d)d\text{ is even},\\
d\neq p-1\text{ if } p\geq3 \text{ and } i=2,\\
d\neq1\text{ if }i=a \text{ and }i\geq3
\end{array}\right\}.
\]
Let $m_i=\max\mathcal W_i$ and $M_i=\sum_{j=1}^{i}m_j$.
Without the two exclusions, the largest elements of the sets $\mathcal W_i$ have sum $S_p(a)$.
When $p\geq3$ is odd, the first exclusion applies only at $i=2$.
The parity condition then replaces $p+1$ by $p-1$, decreasing the sum by two.
For $i=a\geq3$, the second exclusion changes the largest element only when $p=2$ and $a$ is odd.
It then replaces $a+1$ by $a$.
Thus
\[
M_a=\sum_{i=1}^{a}m_i=S_p^\kappa(a)\geq n.
\]

To show that $n$ can be written as $\sum_{i=1}^{a}w_i$ with $w_i\in\mathcal W_i$, we prove by induction on $j$ the stronger statement that every integer in $[0,M_j]$ can be written as $\sum_{i=1}^{j}w_i$ with $w_i\in\mathcal W_i$.
For $j=1$, we have $\mathcal W_1=\{0,1,\ldots,p\}$, so the assertion holds.

Now let $j\geq2$, and suppose that the first $j-1$ sets generate every integer in $[0,M_{j-1}]$.
Write the elements of $\mathcal W_j$ as $x_0=0$ and
\[
j=x_1<x_2<\cdots<x_r=m_j.
\]
The positive elements of $\mathcal W_j$ have consecutive differences at most two.
Indeed, every even value of $d$ is admissible, and the two exclusions either remove the largest element or remove the element corresponding to $d=1$, neither of which creates a gap larger than two.
By the induction hypothesis, choosing $w_j=x_t$ produces every integer in $x_t+[0,M_{j-1}]$.
We need to show that
\[
\bigcup_{t=0}^r\left(x_t+[0,M_{j-1}]\right)=[0,M_j].
\]
It is enough to show that
\begin{equation}\label{eq:vertex_xt_xt+1}
x_{t+1}\leq x_t+M_{j-1}+1
\quad\text{for all }0\leq t<r.
\end{equation}
Since $j-1\in\mathcal W_{j-1}$, we have $M_{j-1}\geq m_{j-1}\geq j-1$.
Thus \eqref{eq:vertex_xt_xt+1} holds for $t=0$.
Moreover, $M_{j-1}\geq m_1=p$.
For $t\geq1$, the consecutive-difference bound gives
\[
x_{t+1}-x_t\leq2\leq p+1\leq M_{j-1}+1.
\]
Hence \eqref{eq:vertex_xt_xt+1} also holds for $t\geq1$.
This completes the induction.

Choose $w_i\in\mathcal W_i$ with $\sum_{i=1}^aw_i=n$.
Whenever $w_i>0$, write $w_i=i+d$ and construct a $d$-regular graph $R_i$ on $i+d$ vertices as follows.
\begin{equation}\label{eq:Ri}
R_i=
\begin{cases}
K_{d+1}, & i=1,\\
\overline{K_i}, & i\geq2,\ d=0,\\
C_{i+d}^{d/2}, & i\geq2,\ d\geq2\text{ even},\\
\overline{C_{i+d}^{(i-1)/2}}, & i\geq3,\ d\text{ odd}.
\end{cases}
\end{equation}
The first condition defining $\mathcal W_i$ makes every case well defined.
We claim that
\begin{equation}\label{blockconnectivitybounds}
	\kappa(R_i)\geq i+d-a,
	\qquad
	\kappa(\overline{R_i})\geq i+d-p.
\end{equation}
We note that if $\kappa(R_i)=d$, the first inequality in \eqref{blockconnectivitybounds} holds since $i\leq a$. 
If $\kappa(\overline{R_i})=i-1$, the second inequality in \eqref{blockconnectivitybounds} holds since $d\leq p-1$.

In the first two cases of \eqref{eq:Ri}, $\kappa(R_i)=d$ and $\kappa(\overline{R_i})=i-1$.
In the third case when $i\geq 3$ and in the fourth case when $d\geq 3$, Lemma~\ref{cyclepowerconnectivity} gives $\kappa(R_i)=d$ and $\kappa(\overline{R_i})=i-1$.

In the third case when $i=2$, the graph is a complete graph minus a perfect matching.
Thus $\kappa(R_i)=d$ and $\kappa(\overline{R_i})=0$.
Here $d\geq2$ implies $p\geq3$.
If $p$ is even, parity gives $d\leq p-2$.
If $p$ is odd, the second restriction on $\mathcal W_i$ excludes $d=p-1$, so again $d\leq p-2$.
Therefore $i+d-p\leq0$, and \eqref{blockconnectivitybounds} holds.

In the fourth case when $d=1$, $R_i$ is a perfect matching.
Hence $\kappa(R_i)=0$ and $\kappa(\overline{R_i})=i-1$.
The third restriction on $\mathcal W_i$ gives $i\leq a-1$.
Thus $i+d-a\leq0$, and \eqref{blockconnectivitybounds} holds.
	
Let $H=\bigvee_{w_i>0}R_i$.
Every vertex of $R_i$ has degree $n-i$ in $H$, so these blocks are the degree classes of $H$ and
\[
H_\delta=\bigvee_{w_i>0}\overline{R_i}.
\]
For a join $J=\bigvee_j J_j$,
\[
\kappa(J)=\min_j\bigl(|V(J)|-|V(J_j)|+\kappa(J_j)\bigr).
\]
The join formula and \eqref{blockconnectivitybounds} give
\[
\kappa(H)\geq n-a,
\qquad
\kappa(H_\delta)\geq n-p.
\]
Therefore $\kappa(H)\kappa(H_\delta)\geq(n-a)(n-p)$.
The upper bound in \eqref{exactKformula}, together with the choice of $p$ and $a$, gives the reverse inequality.
Hence equality holds.
\end{proof}

\begin{example}
Let $n=16$.
For $1\leq p\leq5$, direct evaluation gives
\[
\begin{array}{c|ccccc}
p&1&2&3&4&5\\ \hline
a_{16}^\kappa(p)&6&5&4&4&3\\
(16-p)\bigl(16-a_{16}^\kappa(p)\bigr)&150&154&156&144&143.
\end{array}
\]
For $p\geq6$, the product is at most $(16-6)(16-1)=150$.
Thus the maximum in \eqref{exactKformula} is $156$, attained at $p=3$ and $a=a_{16}^\kappa(3)=4$.
Here
\[
S_3^\kappa(4)=S_3(4)-2=(3+4+5+6)-2=16.
\]
The first exclusion in the definition of $\mathcal W_i$ removes the value $4$ from $\mathcal W_2$, giving
\[
\mathcal W_1=\{0,1,2,3\},
\qquad
\mathcal W_2=\{0,2\},
\qquad
\mathcal W_3=\{0,3,4,5\},
\qquad
\mathcal W_4=\{0,4,6\}.
\]
The largest elements, $3,2,5,6$, sum to $16$.
Choose $(w_1,w_2,w_3,w_4)=(3,2,5,6)$.
The corresponding graph is
\[
H=K_3\vee\overline{K_2}\vee C_5\vee C_6.
\]
The construction gives $\kappa(H)=12$ and $\kappa(H_\delta)=13$, so $H$ attains the product $156$.
\end{example}


\section{Conclusion and Discussion}\label{end}

This paper establishes sharp lower and upper bounds for the sums and products of the minimum degree, maximum degree, vertex connectivity, and edge connectivity of a graph and its $\delta$-complement.
The statements include the necessary small-order qualifications and specify equality cases for every sharp bound.
The degree-class structure controls all four parameters, while complete multipartite graphs and joins of regular graphs provide the extremal constructions.

Further directions include relations for the $\delta'$-complement $G_{\delta'}=\overline{G_\delta}$~\cite{pai}.
It would also be natural to study corresponding relations under subgraph complementation~\cite{kaminski2009} and Seidel switching~\cite{van1966}.
Other Nordhaus--Gaddum questions may consider pairs of different invariants evaluated on $G$ and $G_\delta$.
A broader collection of such relations is surveyed in~\cite{ingrid}.

%
%

\section*{Funding}
Supakorn Srisawat was supported by a Graduate Fellowship (Research Assistant), Faculty of Science, Prince of Songkla University, contract no. 1-2565-02-028.

%
%
%

\bibliography{reference}

\end{document}